\documentstyle[12pt]{article} \newcommand{\nc}{\newcommand}  \nc{\ov}{\over}
\nc{\ra}{\rightarrow} \nc{\iy}{\infty} \nc{\inv}{^{-1}} \nc{\la}{\lambda}
\nc{\be}{\begin{equation}} \nc{\ee}{\end{equation}} \nc{\cd}{\cdots}  \nc{\al}{\alpha}
\renewcommand{\Pr}{{\rm Pr}\,} \nc{\ph}{\varphi} \nc{\tr}{{\rm tr}\,} \nc{\dl}{\delta}
\nc{\noi}{\noindent} \textwidth=6.5in \topmargin=-.5in \textheight=9in
\renewcommand{\sp}{\vspace{1ex}} \nc{\ps}{\psi} \nc{\si}{\sigma}
\begin{document}

\begin{center}{\bf On Convergence of Moments for Random Young Tableaux\\
and a Random Growth Model}\end{center}

\begin{center}{{\bf Harold Widom}\footnote{Supported by National Science 
Foundation grant DMS-9732687.}\\
{\it Department of Mathematics\\
University of California, Santa Cruz, CA 95064\\
e-mail: widom@math.ucsc.edu}}\end{center}

\begin{abstract} In recent work of Baik, Deift and Rains convergence of moments
was established for the limiting joint distribution of the lengths of the first $k$
rows in random Young tableaux. The main difficulty was obtaining a good
estimate for the tail of the distribution and this was accomplished through a
highly nontrival Riemann-Hilbert analysis. Here we give a simpler derivation. The
same method is used to establish convergence of moments for a random growth model.
\end{abstract}

Since the paper of Baik, Deift and Johansson \cite{BDJ1} in which the authors determined
the limiting distribution
of the length of the longest increasing subsequence in a random permutation, or equivalently
the length $\la_1$ of the first row in a random Young tableau, there have been a variety of
extensions and generalizations. In \cite{BDJ2} the same authors determined the limiting 
distribution of the length $\la_2$ of the second row and Borodin, Okounkov and Olshansky 
in \cite{BOO} and Johansson in 
\cite{J2} determined the limit of the joint distribution of $\la_1,\cd,\la_k$ for any
$k$. The result was that for fixed $x_1,\cd,x_k$ the limit
\[\lim_{N\ra\iy}\Pr\left(\la_i\le 2\sqrt N + x_i\, N^{1/6},\,i=1,\cd,k\right)\]
exists and equals the limiting joint distribution of the largest $k$ 
eigenvalues in the Gaussian
unitary ensemble of random matrices. (Here $N$ was the number of boxes in the Young
tableaux, which were given the Plancherel measure.)

In the end it turned out to be the question of the asymptotics of certain Fredholm
determinants, and their derivatives, associated with the operator $K_n$ 
acting on $\ell^2(n,n+1,\cd)$ whose matrix entries are given by
\be K(i,j)=\sum_{k=1}^{\iy}\left({\ph_-\ov\ph_+}\right)_{i+k} \,
\left({\ph_+\ov\ph_-}\right)_{-k-j}.\label{Sn}\ee
Subscripts denote Fourier coefficients and here $\ph_{\pm}(z)=e^{tz^{\pm1}}$, the
Wiener-Hopf factors of $\ph(z)=e^{t(z+z\inv)}$.

In order to prove convergence of moments of the distributions one has to have uniform
estimates for the probabilities and the only difficulties occur at $-\iy$. Convergence 
of moments was established for $\la_1$
and $\la_2$ in \cite{BDJ1,BDJ2} and, for the joint distribution, in \cite{BDR}.
Essentially what one has to show is that each 
\[{d^k\ov dv^k}\det(I-vK_n)\Big|_{v=1}\]
is very small if $2t-n$ is large. In the last cited paper an integral 
operator  was introduced which was an alternative to $K_n$
and which was amenable to a highly nontrivial Riemann-Hilbert analysis which gave the
necessary estimates. Here we shall show that a general determinant inequality allows one to 
derive the estimates more easily.

In \cite{GTW} the authors studied a random growth model called oriented digital boiling
and determined the limiting shape and its fluctuations. In this case the
probability in question, $\Pr(H\le h)$ where $H$ is a certain random variable, is
equal to $\det(I-K_h)$ where now 
\[\ph(z)=(1+z)^{n}\,(1-rz\inv)^{-m},\]
with Wiener-Hopf factors
\be\ph_+(z)=(1+z)^{n},\ \ \ \ph_-(z)=(1-rz\inv)^{-m}.\label{WH}\ee
The asymptotics were determined in various regimes. In the GUE universal regime of 
\cite{GTW} $n=\al m$ with $\al$ fixed and $\al r<1$.
It was shown that for a certain constant $c$, the ``time constant'', the limit
\[\lim_{m\ra \iy}\Pr(H\le cm+s\,m^{1/3})\]
exists and equals the limiting distribution function for the largest eigenvalue in the Gaussian
unitary ensemble. To
establish convergence of moments the main problem is to show that the
determinant is very small, uniformly in $m$, for large negative $s$. (Estimates established 
in \cite{GTW} show that the probability approaches 1 exponentially for large positive $s$,
uniformly in $m$.)

The inequality which allows us to obtain adequate estimates easily is the following.\sp

\noi{\bf Lemma 1}. If $K$ is a trace class operator all of whose eigenvalues are
real and $\le 1$ then
\[\det(I-K)\le e^{-\tr K}.\]

\noi{\bf Proof}. If $\la\le 1$ then $0\le 1-\la\le e^{-\la}$. Take the product over all 
eigenvalues of $K$.\sp

We first apply the lemma to the Young tableaux distributions. As observed in \cite{BDR}, 
since all eigenvaues of $K_n$ in this case are real and lie between 0 and 1,
it is enough to estimate $\det(1-v K_n)$ for any fixed $v\in (0,1)$. If $2t=n+sn^{1/3}$
the authors obtained bounds of the form
\[ \det(I-v K_n)\le\left\{\begin{array}{ll}e^{-\dl s^{3/2}}&{\rm if}\ s\le n^{2/3}\\&\\
e^{-\dl t}&{\rm if}\ s\ge n^{2/3}\end{array}\right.\]
for some $\dl>0$. (These are equivalent to inequalities (5.10) and (5.11) of \cite{BDR}.)
These bounds will follow from the corresponding estimate 
for the trace.\sp

\noi{\bf Proposition 1}. We have for some $\dl>0$
\[\tr K_n\ge\left\{\begin{array}{ll}\dl s^{3/2}&{\rm if}\ s\le n^{2/3}\\&\\
\dl t&{\rm if}\ s\ge n^{2/3}.\end{array}\right.\]\sp

\noi{\bf Proof}. We have 
\[\tr K_n=\sum_{k=1}^{\iy}k\,\left({\ph_-\ov\ph_+}\right)_{n+k}\,
\left({\ph_-\ov\ph_+}\right)_{-n-k}\]
and here $(\ph_-/\ph_+)_k=(\ph_+/\ph_-)_{-k}=J_k(2t)$,
where $J_k$ is the bessel function. Therefore
\be\tr K_n=\sum_{k=1}^{\iy}k\,|J_{n+k}(2t)|^2.\label{trsum}\ee

The uniform asymptotics for the Bessel function (\cite{E}, sec. 4.8) give for $u>1$
\[J_k(ku)=\left({1\ov 2}k^{2/3}u\psi'(u)\right)^{-1/2}{\rm Ai}(-k^{2/3}\psi(u))\,
(1+O(k\inv)),\]
where ${2\ov3}\psi(u)^{3/2}=\int_1^u(1-t^{-2})^{1/2}\,dt$.
(The error term has to be modified near the zeros of the Airy function, that is, 
when ${2\ov 3}k\psi(u)^{3/2}$ is near the set $-{\pi\ov4}+\pi{\bf Z}$.) 
It is easy to see that if
\[2t=k+s\,k^{1/3}\]
then $k^{2/3}\psi(2t/k)$ tends to $+\iy$ as $k,\,s\ra+\iy$. Hence, from the asymptotics of
the Airy function we have for large $k$ and $s$ and some constant $\dl>0$
\be|J_k(2t)|\ge \dl\,k^{-1/2}\left({2t\ov k}\psi'\left({2t\ov k}\right)\right)^{-1/2}
\psi\left({2t\ov k}\right)^{-1/4}\label{bound}\ee
as long as ${2\ov 3}k\psi(2t/k)^{3/2}$ is bounded away from $-{\pi\ov4}+\pi{\bf Z}$.

Now $t$ and $n$ are related by
\[2t=n+s\,n^{1/3}.\] 
We replace $k$ by $n+k$ in (\ref{bound}) and choose those $k$ such that
\be \al_1\,s\,n^{1/3}\le k\le \al_2\,s\,n^{1/3},\label{k}\ee
where $\al_1<\al_2<1$. Then we find that 
\[{2t\ov n+k}-1\approx\left\{\begin{array}{ll}sn^{-2/3}&{\rm if}\ s\le n^{2/3}\\&\\
1&{\rm if}\ s\ge n^{2/3},\end{array}\right.\]
where ``$\approx$'' means ``is bounded above and below by positive constants times''. 
{}From this and (\ref{bound}) we deduce that with another $\dl>0$\ \footnote{
We should be careful
because of the requirement that ${2\ov 3}(n+k)\psi(2t/(n+k))^{3/2}$ be bounded away from 
$-{\pi\ov4}+\pi{\bf Z}$. If $\al_1$ in (\ref{k}) is very close to
1 then two consecutive values of $(n+k)\psi(2t/(n+k))^{3/2}$ 
differ by at most a small constant. Thus with such a choice of $\al_1$ we can be
assured that the estimates hold for a large fraction of the $k$ satisfying (\ref{k}).}
\[|J_{n+k}(2t)|\ge\left\{\begin{array}{ll}\dl s^{-3/4}&{\rm if}\ s\le n^{2/3}\\&\\
\dl n^{-1/6}s^{-1/2}&{\rm if}\ s\ge n^{2/3}.
\end{array}\right.\]
And from this, summing in (\ref{trsum}) over only those $k$ satisfying (\ref{k}), we deduce 
that for $s\le n^{2/3}$ the trace of 
$K_n$ is at least a constant times $s^{1/2}n^{2/3}$
which is in turn at least a constant times $s^{3/2}$ and for $s\ge n^{2/3}$
it is at least a constant times $sn^{1/3}$ which is in turn at least a constant 
times $t$. This proves the proposition.\sp

We turn now to the growth model where $K_h$ is defined in terms of the Wiener-Hopf
factors (\ref{WH}).\footnote{Now the operator is not trace class on the
usual $\ell^2$ space. From
the integral formula for the matrix entries given at the beginning of the proof of
Lemma~3 one sees that they grow in $i$ as $i^{n-1}$ and decrease expontially in $j$ as
$j\ra\iy$. So we must take our $\ell^2$ space to be weighted. If the weight
function is $w(i)=2^{-i}$, for example, the operators $K_h$ are trace class.}\sp

\noi{\bf Lemma 2}. The eigenvalues of $K_h$ are all real and less than 1.\sp

\noi{\bf Proof}. An expression for $K(i,j)$ when
\[\ph(z)=(1+\xi z)^{n}\,(1+\xi z\inv)^{-m}\]
is given by formula (4.9) of \cite{BO}.\footnote{Our $K(i,j)$ differs from that of
\cite{BO} by a factor $(-1)^{i+j}$. This does not affect the eigenvalues and so we
make believe the factor is not there.} The $i,j$ entry 
of our operator is equal $(i\,r^{1/2})^{i-j}$
times the $i,j$ entry of this one in which $\xi=i\,r^{1/2}$. From this and the 
formula given there we see that our $K(i,j)$ is of the form $a_i\,L(i,j)\,b_j$ where
$a_i,\,b_j\ge0$ and $L(i,j)$ is symmetric.
Hence all eigenvalues of $K_h$ are real. Next,
notice that $\det\,(I-K_h)$ is always positive,
being a nonzero probability. So no eigenvalue of $K_h$ can equal 1. 
The eigenvalues
being continuous functions of $r$, and the eigenvalues being all 0 when $r=0$, we
see that no eigenvalue can be $\ge1$.\sp

Recall that $n=\al m$ with $\al$ fixed, $\al r<1$ and 
\[\Pr(H\le cm+s\,m^{1/3})=\det\,(I-K_h)\]
where $h=cm+s\,m^{1/3}$. The time constant $c$ is given by the formula
\[c={1\ov 1+r}((1-\al)r+2\sqrt{\al r}).\]
Since the only problem occurs when $s$ is large negative  we change $s$ to $-s$ and 
think of $s$ as large and positive. So we have to estimate $\det\,(I-K_h)$
when $h=cm-s\,m^{1/3}$. By Lemmas~1 and 2 it suffices to find a lower bound for the trace.\sp

\noi{\bf Lemma 3}. If $s$ is sufficiently large and $sm^{-2/3}$ sufficiently small there 
exists a constant
$\delta>0$ such that $\tr K_h\ge\delta s^{3/2}$.\sp

\noi{\bf Proof}. Given functions $f$ and $g$ we have
\[\sum_{k=1}^{\iy}f_{i+k} \,g_{-k-j}=-{1\ov4\pi^2}\int\int f(z_1)\,g(z_2)\,
\sum_{k=1}^{\iy}z_1^{-i-k-1}\,z_2^{k+j-1}\,dz_1\,dz_2\]
\[=-{1\ov4\pi^2}\int\int f(z_1)\,g(z_2)\,
{z_1^{-i-2}\,z_2^{j}\ov 1-z_2/z_1}\,dz_1\,dz_2,\]
where the integrals are taken over circles about the origin on which $|z_2|<|z_1|$.
If we consider these as matrix elements of an operator acting on
$\ell^2(h,\,h+1,\cdots)$ then the trace of the operator equals
\[-{1\ov4\pi^2}\int\int f(z_1)\,g(z_2)\,
{z_1^{-h}\,z_2^{h}\ov (z_1-z_2)^2}\,dz_1\,dz_2.\]
In our case we set
\[f(z)=(1+z)^{-n}\,(1-rz\inv)^{-m},\ \ \ g(z)=(1+z)^{n}\,(1-rz\inv)^{m}\]
and define, following the notation of \cite{GTW},
\[\ps(z,c')=(1+z)^{n}\,(1-rz\inv)^{m}\,z^{c'm}.\]
If $h=c'm$ we have
\[\tr K_h=-{1\ov4\pi^2}\int\int {\ps(z_2,c')\ov \ps(z_1,c')} {dz_1\,dz_2\ov (z_1-z_2)^2}.\]
Here we should have $|z_1|<1$ and $|z_2|>r$, so $r<1$. The formulas for the matrix entries
and trace hold for all $r$
by analytic continuation if the $z_2$ contour is inside the $z_1$ contour and the latter 
has $r$ on the inside and 
$-1$ on the outside; the $z_2$ contour has 0 on the inside.

Now $h=cm-sm^{1/3}$, so $c'=c-sm^{-2/3}$,
and we want a lower bound for the trace. We use steepest descent
on each of the two integrals as if the factor $(z_1-z_2)^{-2}$ were not there. 
When $c'=c$ the function
\[\si(z,c')={1\ov m}\log\,\ps(z,c')=\al\,\log\,(1+z)+\log\,(r-z)+(c'-1)\log\,(-z)\]
has vanishing first and second derivatives at the single critical point $u_c\in (-1,0)$ and
the two steepest descent curves are as shown in Fig. 6 of \cite{GTW}. The one for the
$z_2$ integration, denoted here by $C_2$, is closed and contains 0 on the inside while $C_1$,
the one for the $z_1$
integration, has two arms going to infinity and separates 0 from $-1$. Our original
contours can be deformed to these.

For $c'\ne c$ there are in general two critical points $u_{c'}^{\pm}$ given by
\be u^{\pm}_{c'}={-(1-r)c'-(1-\al)r \pm
\sqrt{\left((1+r)\,c'+(\al-1)\,r\right)^2-4\al\, r}\ov2(\al+c')}.\label{u}\ee
The square root is zero when $c'=c$.
When $c'$ decreases from $c$ the square roots are, at least for a while, purely 
imaginary.\footnote{There are cases when the square roots are real for some $c'<c$, and the 
following descussion would have to be modified when this happens. See Remark~1 following 
the proof of the lemma. This is the reason for our assumption that
$sm^{-2/3}=c-c'$ is sufficiently small.} Denote by
$u_{c'}^+$ the one in the upper half-plane and by $u_{c'}^-$ the one in the lower 
half-plane. As we can see the steepest descent contours $C_1$ and $C_2$ now cross
and so $C_2$ is no longer inside inside $C_1$. (Also, the resulting double integral
is not absolutely convergent but must be interpreted as a principal value). Therefore
our original double integral is not equal to the integral over $C_1\times C_2$. In fact,
we can see that our original double integral (without its factor $-1/4\pi^2$) is equal to 
the integral over $C_1\times C_2$ plus
\be\int_{u_{c'}^-}^{u_{c'}^+}\int_{\Gamma} {\ps(z_2,c')\ov \ps(z_1,c')} 
{dz_1\,dz_2\ov (z_1-z_2)^2},\label{psiints}\ee
where the $z_2$ path of integration is to the right of $-1$ and to the left of 0 and
the loop $\Gamma$ enclose all point of the $z_2$ path. The inner integral is equal to
\[-2\pi i{\ps'(z_2,c')\ov\ps(z_2,c')}\]
and integrating this and dividing by the factor $-4\pi^2$ gives
\[-{1\ov 2\pi i}[\log\ps(u_{c'}^+,c')-\log\ps(u_{c'}^+,c')]=
-{m\ov 2\pi i}[\si(u_{c'}^+,c')-\si(u_{c'}^-,c')],\]
where in the computation of the difference of logarithms we go over a path 
to the right of $-1$ and to the left of 0.

It follows from the displayed formula after (3.17) in \cite{GTW} that
\[\si(u_{c'}^+,c')-\si(u_{c'}^-,c')=\int_c^{c'} \log {u_{\gamma}^+\ov 
u_{\gamma}^-}\,d\gamma.\]
The logarithm in the integral is positive purely imaginary and $c'<c$, so
the left side is negative purely imaginary. Using (\ref{u}) we see that in fact 
this expression
times its factor $-m/2\pi i$ is at least a positive constant times $m(c-c')^{3/2}=s^{3/2}$.

It remains to estimate
\[\int_{C_2}\int_{C_1} {\ps(z_2,c')\ov \ps(z_1,c')} 
{dz_1\,dz_2\ov (z_1-z_2)^2},\]
and we shall show that it is $O(1)$. It can be seen that $C_1$ closes at $-1$ while
$C_2$ closes at 0. Denote by $C_i^+$ the portions of these contours in the
upper half-plane and  by $C_i^-$ the portions in the lower half-plane. The distances
between $C_1^{\pm}$ and $C_2^{\mp}$ are of order $(c-c')^{1/2}$.
Lemma~6 of \cite{GTW} tells us that $\si''(u_{c'}^{\pm},c')$ is of
the order $\pm(c-c')^{1/2}$. It follows that the integrals over $C_1^{\pm}\times C_2^{\mp}$
are of order $(c-c')^{-1}m^{-1}(c-c')^{-1/2}=s^{-3/2}$.
So these integrals are $O(1)$ for $s\ge1$ and we may restrict attention to the 
integrals over $C_1^{\pm}\times C_2^{\pm}$.

Considering the integral over $C_1^{+}\times C_2^{+}$ we can confine attention to an
arbitrarily small neighborhood of $(u_{c'}^+,u_{c'}^+)$. (The integral over the complement
of any such neighborhood is exponentially small.) For notational convenience
we write $u$ for $u_{c'}^+$. In our neighborhood we have
\[{\ps(z_2,c')\ov \ps(z_1,c')}=e^{m[\si''(u,c')((z_2-u)^2-(z_1-u)^2)+
O(|z_2-u|^3+|z_1-u|^3)]}\]
\[=e^{m\si''(u,c')((z_2-u)^2-(z_1-u)^2)}\left(1+O(m(|z_2-u|^3+|z_1-u|^3))\right).\]
Because our contours cross at a positive angle we have
$(z_i-u)^3/(z_1-z_2)^2=O(|z_i-u|)$ and it follows that the contribution of the $O$ term
on the right is at most a constant times
\[{m\ov (m^{1/2}(c-c')^{1/4})^3}=s^{-3/4}.\]
So this contribution is also $O(1)$ for $s\ge1$.

That leaves the main integral
\[\int_{C_2^+}\int_{C_1^+} {e^{m\si''(u,c')((z_2-u)^2-(z_1-u)^2)}\ov (z_1-z_2)^2}
dz_1\,dz_2\]
which, recall, is interpreted as a principal value. Recall also that we confine attention
to a neighborhood of $(u,u)$. If we make the variable
changes $z_i-u=\zeta_i/(m\si''(u,c'))^{1/2}$ with an appropriate choice of square root
the above becomes
\[\int_{\Gamma_2}\int_{\Gamma_1}{e^{(\zeta_2^2-\zeta_1^2)}\ov (\zeta_1-\zeta_2)^2}
d\zeta_1\,d\zeta_2,\]
where $\Gamma_i$ are long contours crossing at $(0,0)$, with $\Gamma_1$
horizontal there and $\Gamma_2$ vertical. The integral outside a neighborhood of $(0,0)$
is bounded. In a neighborhood of $(0,0)$ we can replace the exponential by 1 with error
$O(1)$ and the principal value integral resulting after this replacement is $O(1)$.
This establishes the lemma. \sp

\noi{\bf Remark 1}. We know that $c$ is a value of $c'$ such that there 
is only one critical point. But if $2\sqrt{\al r}\le (1-\al)r$ there
will be another nonnegative value of $c'$ where there is only one critical point, namely
\[{1\ov 1+r}((1-\al)r-2\sqrt{\al r}).\]
And for smaller $c'$ there will be two real critical points. So the nature of the
curves $C_i$ will be different. In particular for $c'$ less than this value $C_2$ is
outside $C_1$. What happens in these cases is that instead of integrating from $u_{c'}^-$
to $u_{c'}^+$ in (\ref{psiints}) the outer integral is a closed loop with the result that the integral with
its factor $-1/4\pi^2$ is equal to $n$ exactly. However, concerning ourselves with this 
is unnecessary since we know enough already.\sp

\noi{\bf Proposition 2}. There exists a $\dl>0$ such that $\det(I-K_h)\le e^{-\dl s^{3/2}}$
for sufficiently large $s$.\sp

\noi{\bf Proof}. For sufficiently small $sm^{-2/3}$, say 
$sm^{-2/3}\le \eta$, Lemmas~1--3 give the 
stated estimate. In particular for $sm^{-2/3}=\eta$ the determinant is at most
$e^{-\delta s^{3/2}}=e^{-\delta\eta^{3/2} m}$. Since the determinant is a nonincreasing
function of $s$ this bound holds also for $sm^{-2/3}>\eta$. But since $h=cm-sm^{1/3}\ge0$
we have $m\ge c^{-3/2}s^{3/2}$ and therefore the stated bound holds also for these $s$, 
with $\delta$ replaced by $\delta(\eta/c)^{3/2}$.\sp

\noi{\bf Remark 2}. Baik, Deift, McLaughlin, Miller and Zhou \cite{BDMMZ} have recently
proved convergence of moments, using Riemann-Hilbert techniques, for a
growth model of Johansson \cite{J1}. There an anologous determinant arises where
$K(i,j)$ is given by (\ref{Sn}) with
\[\ph_+(z)=(1+t\,z)^M,\ \ \ \ph_-(z)=(1+t\,z\inv)^N.\]
We thank the authors for alerting us to their results.
Very likely their methods could also be used for our growth model. 

\noi{\bf Remark 3}. The proof of Lemma 2 is not very satisfactory. There should be a reason for
its truth and it should be a special case of a more general fact. There may
be some merit in the following\sp

\noi{\bf Conjecture}. Let $r_i$ and $s_j$ be finitely many nonegative real numbers.
Then the eigenvalues of the operator $K_n$ acting on $\ell^2(n,n+1,\cd)$ whose matrix 
entries are  given by (\ref{Sn}) with
\[\ph_+(z)=\prod_i \,(1+r_i\,z),\ \ \ \ph_-(z)=\prod_j\,(1-s_j\,z\inv)\inv\]
are real and lie between 0 and 1.
\sp 

Numerical evidence supports the truth of the conjecture.\sp\sp

\end{document}